\documentclass[11pt]{article}

\usepackage{latexsym,color,graphicx,amsmath,amssymb,amsthm,graphics}
\usepackage{subfig}

\def\reals{{\mathbb R}}
\def\R{{\mathbb R}}

\def\N{{\mathbb N}}

\newtheorem{theorem}{Theorem}%[section]
\newtheorem{lemma}[theorem]{Lemma}

\theoremstyle{definition}
\newtheorem{definition}[theorem]{Definition}

\def\N{\mathbb {N}}

\def\E{\mathsf{E}}

\def\a{\alpha}

\def\({\big (}
\def\){\big )}

\def\le{\leqslant}
\def\ge{\geqslant}
\def\_phi{\varphi}

\def\R{{\mathbb R}}

\begin{document}

\title{On the number of unit-area triangles spanned by convex grids in the plane\thanks{%
Work on this paper by Orit E. Raz and Micha Sharir was supported by Grant 892/13 from the Israel Science Foundation,
and by the Israeli Centers of Research Excellence (I-CORE) program (Center No.~4/11).
Work by Micha Sharir was also supported by Grant 2012/229 from the U.S.--Israel Binational Science Foundation,
and by the Hermann Minkowski-MINERVA Center for Geometry at Tel Aviv University.
Work by Ilya D. Shkredov was supported by Russian Scientific Foundation RSF 14--11--00433.
}}

\author{
Orit E. Raz\thanks{%
School of Computer Science, Tel Aviv University,
Tel Aviv 69978, Israel.
{\sl oritraz@post.tau.ac.il} }
\and
Micha Sharir\thanks{%
School of Computer Science, Tel Aviv University,
Tel Aviv 69978, Israel.
{\sl michas@post.tau.ac.il} }
\and
Ilya D. Shkredov\thanks{%
Division of Algebra and Number Theory, Steklov Mathematical Institute, ul. Gubkina 8, Moscow, 119991, Russia
and IITP RAS,  Bolshoy Karetny per. 19, Moscow, Russia, 12799.
{\sl ilya.shkredov@gmail.com}}
}
%\date{}
\maketitle

\begin{abstract}
A finite set of real numbers is called convex if the differences between consecutive
elements form a strictly increasing sequence. We show that, for any pair of convex
sets $A, B\subset\reals$, each of size $n^{1/2}$, the convex grid $A\times B$ spans
at most $O(n^{37/17}\log^{2/17}n)$ unit-area triangles. This improves the best known
upper bound $O(n^{31/14})$ recently
obtained in \cite{RS}. Our analysis also applies to more general families of sets
$A$, $B$, known as sets of Szemer\'edi--Trotter type.
\end{abstract}

%%%%%%%%%%%%%%%%%%%%%%%%%%%%%%%%%%%%%%%%%%%%%%%%%%%%%%%%%%%%%%%%%%%

\section{Introduction}

The problem considered in this paper is to obtain
a sharp upper bound on the maximum number of unit-area triangles spanned by $n$ points in the plane.
The problem has been posed by Oppenheim in 1967 (see \cite{EP95}) and has been studied since then
in a series of papers \cite{AS10,DST08,EP71,PS92,RS}.
The currently best known upper bound, for an arbitrary set of $n$ 
distinct points in the plane, is due to the first two authors~\cite{RS}:
\begin{theorem}[{\bf Raz and Sharir}~{\bf \cite{RS}}]\label{main1}
The number of unit-area triangles spanned by $n$ points in the plane is $O(n^{20/9})$.
\end{theorem}

In this paper we study a variation of the problem, first considered in~\cite{RS},
concerning the number of unit-area triangles spanned by points in a {\em convex grid}.
That is, the input set is of the form $A\times B$, where $A$ and $B$
are convex sets  of $n^{1/2}$ real numbers each;
a set $X=\{x_1,\ldots,x_n\}$, with $x_1<x_2<\cdots<x_n$, of real numbers is called \emph{convex} if
$$
x_{i+1}-x_i>x_i-x_{i-1},
$$
for every $i=2,\ldots,n-1$. See \cite{ENR, Gar00, GaK05, Heg86, IKRT06, Li11, Sch, SS11, Sol09} for more details and properties  of convex sets.

We establish the following improvement of Theorem~\ref{main1} for convex grids.
%%%%%%%%%%%%%%%%%%%%%%%%%%%%%%%%%%%%%%%%%%%%%%%%%%%%%%%
\begin{theorem}\label{main3}
Let $S=A\times B$, where $A,B\subset \R$ are convex sets of size $n^{1/2}$ each.
Then the number of unit-area triangles spanned by the points of $S$ is $O(n^{37/17} \log^{2/17} n)$.
\end{theorem}
%%%%%%%%%%%%%%%%%%%%%%%%%%%%%%%%%%%%%%%%%%%%%%%%%%%%%%%

As noted, the problem was first considered in the conference paper~\cite{RS}, where the weaker bound $O(n^{31/14})$
has been obtained. The analysis presented here uses the general high-level approach as in \cite{RS}, but develops 
and adapts refined techniques that lead to the improved bound asserted above.

Our result is a somewhat unusual application of the properties of convex sets (and of their generalization, so-called \emph{sets of Szemer\'edi--Trotter type},
defined and discussed below) in a two-dimensional ``purely geometric'' context,
in contrast with their more standard applications in arithmetic combinatorics (such as in \cite{Li11,SS11}). Nevertheless, to
obtain the improved bounds, we need to adapt and develop new properties of convex sets (and of sets of Szemer\'edi--Trotter type)
within their standard arithmetic context.
These properties are presented in Section 2, and culminate in Lemma 9, a technical ``multi-dimensional'' property of such sets,
which we believe to be of independent interest.

The main technical tool used in our analysis is a result of Schoen and Shkredov~\cite{SS11} on difference sets involving convex sets.
For two given finite subsets $X,Y\subset \R$, and for any $s\in\R$, denote by $\delta_{X,Y}(s)$
the number of representations of $s$ in the form $x-y$, with $x\in X$, $y\in Y$.
%%%%%%%%%%%%%%%%%%%%%%%%%%%%%%%%%%%%%%%%%%%%%%%%%%%%%%%
\begin{lemma}[{\bf Schoen and Shkredov~\cite[Lemma 2.6]{SS11}}]\label{lem:conv} %Lemma 7 in [SS11]
Let $X,Y$ be two finite subsets of $\R$, with $X$ convex. Then, for any $\tau\ge 1$, we have
$$
\bigl|\bigl\{s\in X-Y~\mid~\delta_{X,Y}(s)\ge\tau\big\}\big|=O\left(\frac{|X||Y|^2}{\tau^3} \right).
$$
\end{lemma}
%%%%%%%%%%%%%%%%%%%%%%%%%%%%%%%%%% %%%%%%%%%%%%%%%%%%%%%

As a matter of fact, our analysis applies to the more general setup of grids $A\times B$, 
where the sets $A,B$ are of {\it Szemer\'edi--Trotter type}.
We recall their definition from \cite{s_sumsets} (with a slightly different notation, and for the special case of the additive group $\R$).\footnote{Note 
that the factor $d(X)|X|$ in the right-hand side of (\ref{f:SzT-type}) can be replaced 
simply by a function of $X$, $c(X)$, as is done in \cite{s_sumsets}. 
We choose to keep the quantity $|X|$ in the notation, though, to have a better comparison 
between the convex case and the more general case of sets of SzT-type.}
We note that sets of Szemer\'edi-Trotter type have been instrumental in recent works in additive combinatorics, including an improved bound on the
sum-product problem~\cite{KS}; see also \cite{s_balogq}.

%%%%%%%%%%%%%%%%%%%%%%%%%%%%%%%%%%%%%%%%%%%%%%%%%%
\begin{definition}[{\bf Shkredov \cite{s_sumsets}}]
A finite set $X\subset \R$ is said to be of {\it Szemer\'{e}di--Trotter type}
(abbreviated as {\it SzT--type}), with a parameter $\a \ge 1$, if there exists a parameter $d(X)>0$ such that
\begin{equation}\label{f:SzT-type}
\bigl|\bigl\{ s\in X-Y ~\mid~ \delta_{X,Y}(s) \ge \tau \big\}\big| \le\frac{d(X) |X| |Y|^\a}{\tau^{3}}\,,
\end{equation}
 for every finite set $Y\subset \R$ and every real number $\tau \ge 1$.
\label{def:SzT-type}
\end{definition}
%%%%%%%%%%%%%%%%%%%%%%%%%%%%%%%%%%%%%%%%%%%%%%%%%%
In view of Lemma~\ref{lem:conv}, any convex set $X$ is of SzT-type, with the parameter
$\alpha=2$ and with $d(X) = O(1)$; see \cite{SS11} for more details.
 
We note also that any set $X$ is of SzT-type with parameter $\alpha=2$ and $d(X)=|X|$. 
Indeed, it is sufficient to check the condition (\ref{f:SzT-type}) just for $\tau \le \min\{ |X|, |Y|\}\le \sqrt{|X|\cdot |Y|}$. 
Since the left-hand side of \eqref{f:SzT-type} does not exceed $|X| |Y| /\tau$, the claim follows. Thus, in what follows we may assume that $d(X) \le |X|$.

Definition~\ref{def:SzT-type} is rather general, and applies, with nontrivial bounds for $d(X)$, to several other families of sets,
where, for many instances, the condition \eqref{f:SzT-type} is obtained using the Szemer\'edi--Trotter
theorem \cite{ST83} (which is why they are named this way), that we recall next.
%%%%%%%%%%%%%%%%%%%%%%%%%%%%%%%%%%%%%%%%%%%%%%%%%%
\begin{theorem}[{\bf Szemer\'edi and Trotter~\cite{ST83}}]\label{ST}
(i) The number of incidences between $M$ distinct points and $N$ distinct lines in the plane is
$O(M^{2/3}N^{2/3}+M+N)$.
(ii) Given $M$ distinct points in the plane and a parameter $k\le M$, the number of lines incident to at least
$k$ of the points is $O(M^2/k^3+M/k)$. Both bounds are tight in the worst case.
\end{theorem}
%%%%%%%%%%%%%%%%%%%%%%%%%%%%%%%%%%%%%%%%%%%%%%%%%%

Our analysis will actually establish the following generalization of Theorem~\ref{main3}.
%%%%%%%%%%%%%%%%%%%%%%%%%%%%%%%%%%%%%%%%%%%%%%%%%%%%%%
\begin{theorem}\label{SzT2}
Let $S=A\times B$, where $A,B\subset \R$ are of SzT-type with parameter $\alpha=2$, and are of size $n^{1/2}$ each.
Assume further that 
$$
d(A)\le d(B) \le n^{6/47} d(A)^{42/47} (\log n)^{-13/47}.
$$
Then the number of unit-area triangles spanned by the points of $S$ is 
$$O\left(n^{37/17}(d^2 (A)d(B)\log n)^{2/17}\right).
$$
\end{theorem}
%%%%%%%%%%%%%%%%%%%%%%%%%%%%%%%%%%%%%%%%%%%%%%%%%%%%%%
To appreciate the broader scope of Theorem~\ref{SzT2}, let us present additional
examples of sets $X$ of SzT-type with parameter $\a=2$, and with $d(X)\ll |X|$.\\

\noindent{{\bf Examples.}}

\noindent
(i) Let $f$ be a strictly convex or concave function, and let $A\subset \R$ be finite.
Then $f(A)$ is of SzT-type with parameter $\alpha=2$ and with $d(f(A)) = q(A)$, where
\begin{equation}\label{f:q(A)}
    q (A) := \min_{C \neq \emptyset}
        \frac{|A+C|^2}{|A||C|}
        \,.
\end{equation}
By the same token, any set $A$ is of SzT--type with $d(A) = q(f(A))$;
see \cite[Lemma 7]{RRS} and also \cite{Li2,SS11} for more details.

\noindent
(ii) Let $A\subset \R$ be finite, and assume that $|AA| \le M|A|$.
Then $A$ is of SzT-type with $d(A) \le M^2$.
Indeed, assume first that $A\subset \R^+$, and consider the convex function $f(x) = 2^x$ and the set $A':=\log A$.  Then $A=f(A')$. 
By Example (i), $A$ is of SzT-type, with parameter $\alpha=2$ and 
$d(A)=q(A')$. 
We have
$$
d(A)=\min_{C \neq \emptyset} \frac{|A'+C|^2}{|A'||C|}\le \frac{|A'+A'|^2}{|A'||A'|}=\frac{|\log AA|^2}{|A|^2}=\frac{|AA|^2}{|A|^2}\le M^2.
$$
The case in which $A$ also involves negative elements is argued similarly.

\noindent
(iii) Let $A \subset \R^{+}$. %, and $a \in \R \setminus \{0 \}$ be any number fixed.
Then $\log A$ is of SzT-type with parameter $\alpha=2$ and with $d(\log A) = \tilde q(A)$, where
$$
    \tilde q (A) := \min_{a\in\R\setminus\{0\},\;C \neq \emptyset}
        \frac{|(A+a)C|^2}{|A||C|}
        \,;
$$
see \cite{J_RN,B_RN_S}.

In the next section we introduce some notation and preliminary background, and establish the
properties of SzT-type sets which are needed for our analysis.
The proof of Theorem~\ref{SzT2} is then given in Section~\ref{sec:prf}.

%%%%%%%%%%%%%%%%%%%%%%%%%%%%%%%%%%%%%%%%%%%%%%%%%%%%%%%%%%%%%%%

\section{Properties of sets of Szemer\'edi--Trotter type}\label{sec:SzT}

Let $A \subseteq \R$ be finite. By a slight abuse of notation, we use the same letter to
denote the characteristic function $A : \R\to \{0,1\}$ of $A$, with $A(x)=1$ if and only if $x\in A$.
For $k\ge 2$ and a $(k-1)$-tuple $\vec x=(x_1,\dots, x_{k-1}) \in \R^{k-1}$, we put
$$
A_{\vec x}:= A \cap (A-x_1)\dots \cap (A-x_{k-1}).
$$
(For $k=2$ we use the shorter notation $A_x$ for $A\cap (A-x)$, for $x\in\reals$.)
Its characteristic function is then
$$
A_{\vec x}(z)= A(z) A(z+x_1) \cdots A(z+x_{k-1}).
$$
The {\it additive energy} of $A$ (see, e.g., \cite{TV}) is defined as
$$
\E(A):=\sum_{x\in\R}|A_x|^2.
$$
Equivalently,
$$
\E(A)=\Bigl|\{(a,a',b,b')\in A^4\mid a-b=a'-b'\}\Bigl|.
$$
Indeed,
\begin{align*}
\E(A) & = \sum_{x\in\R} \Bigl| \{(a,b)\in A^2 \mid a-b=x\}\Bigr|^2 
 = \sum_{x\in\R} \Bigl| \{a\in A \mid a-x \in A\}\Bigr|^2 
 = \sum_{x\in\R} |A_x|^2 .
\end{align*}

This notion can be extended to yield higher-order variants of the energy.
Specifically, for $k\ge 2$,
the \emph{$k$-order energy} of $A$ is defined as
\begin{equation}\label{f:E_k_preliminalies_B}
\E_k(A):=\sum_{x\in\R}|A_x|^k.
\end{equation}
Thus $\E(A) = \E_2(A)$. Equivalently, as shown in \cite{SV}, we have 
\begin{equation}\label{f:E_k_preliminalies}
\E_k (A)= \sum_{\vec{x}=(x_1,\dots,x_{k-1}) \in \R^{k-1}}  |A_{\vec x}|^2 \,.
\end{equation}
Finally, for $k,l\ge 2$, define
$$
\E_{k,l}(A):=\sum_{\vec{x}=(x_1,\dots,x_{k-1}) \in \R^{k-1}}  |A_{\vec x}|^l \,.
$$
In view of \eqref{f:E_k_preliminalies}, we have $\E_{2,k}=\E_{k,2}(=\E_k)$.
More generally, we have the following property.
\begin{lemma}[{\bf Shkredov and Vyugin~\cite{SV}}]
Let $A\subseteq \R$ be finite, and let $k,l\ge 2$.  Then
$$
\E_{k,l} (A) = \E_{l,k} (A) \,.
$$
\end{lemma}
%%%%%%%%%%%%%%%%%%%%%%%%%%%%%%%%%%%%%%%%%%%%%%%

For a set $A$ of SzT-type (with $\alpha=2$) we have the following simple bound on $\E_3(A)$.
%%%%%%%%%%%%%%%%%%%%%%%%%%%%%%%%%%%%%%%%%%%%%%
\begin{lemma}\label{E3}
Let $A$ be a set of SzT-type with parameter $\a=2$. Then
$$
\E_3(A)= O\left(d(A)|A|^3\log|A|\right).
$$
\end{lemma}
%%%%%%%%%%%%%%%%%%%%%%%%%%%%%%%%%%%%%%%%%%%%%%
\noindent{\bf Proof.}
By definition,
\begin{align*}
\E_3(A)
&=\sum_{s\in\R}|A_s|^3
=\sum_{s\in A-A}\delta_{A,A}^3(s)\\
&=\sum_{\tau=1}^{|A-A|}\tau^3\cdot \big|\{s\mid \delta_{A,A}(s)=\tau\}\big|\\
&=\sum_{\tau=1}^{|A-A|}\tau^3\cdot \Big(\bigl|\{s\mid \delta_{A,A}(s)\ge\tau\}\big|-\bigl|\{s\mid \delta_{A,A}(s)\ge\tau+1\}\big|\Big)\\
&=|A-A|+O\left(\sum_{\tau=2}^{|A-A|}\tau^2\cdot \big|\{s\mid \delta_{A,A}(s)\ge\tau\}\big| \right).
\end{align*}
Using \eqref{f:SzT-type}, this yields
\begin{align*}
\E_3(A) & = |A-A|+ O\left(\sum_{\tau=2}^{|A-A|}\tau^2\cdot \frac{d(A) |A|^3}{\tau^{3}} \right)\\
& = O\left(|A|^2+d(A)|A|^3\log|A-A|\right) =O\left(d(A)|A|^3\log|A|\right),
\end{align*}
as asserted.
$\Box$

We also need the following property.
%%%%%%%%%%%%%%%%%%%%%%%%%%%%%%%%%%%
\begin{lemma}\label{mtaurich}
Let $A$ be a set of SzT-type with parameter $\a=2$, and let $m\ge 3$.
Then
\begin{equation}\label{f:E_kl_C}
    \Big|\Bigl\{ (a_1,\dots, a_m)\in A^m ~\mid~ |(A-a_1) \cap \dots \cap (A-a_m)| \ge \tau \Big\}\Big|
       \le
            \frac{d(A) |A| \E_{m-1,3} (A)}{\tau^3} ,
\end{equation}
\label{l:E_kl_SzT}
 for every $\tau \ge 1$.
\end{lemma}
%%%%%%%%%%%%%%%%%%%%%%%%%%%%%%%%%%%
\noindent{\bf Proof.}
Put
$$
\sigma := \Bigl|\Bigl\{ (a_1,\ldots,a_{m}) \in A^{m} ~\mid~ |(A-a_1) \cap  \cdots \cap (A-a_{m})| \ge \tau \Bigr\} \Bigr| \,.
$$
We have
\begin{align*}
|(A-a_1)  &\cap \dots \cap (A-a_{m})|=\sum_{z\in\reals} A(z+a_1) \cdots A(z+a_m)\\
&= \sum_{z\in\reals} A(z) A(z+a_2-a_1)A(z+a_3-a_1) \cdots A(z+a_m-a_1)\\
&= \sum_{z\in\reals} A(z) A(z+a_2-a_1)A(z+a_2-a_1+a_3-a_2) \cdots A(z+a_2-a_1+a_m-a_2)\\
&= \sum_{z\in\reals} A(z)A_{(a_3-a_2,\ldots,a_m-a_2)}(z+a_2-a_1)\\
&= \delta_{A,A_{(a_3-a_2,\ldots,a_m-a_2)}}(a_1-a_2)\, .
\end{align*}
Hence,
$$
\sigma = \Bigl|\Bigl\{ (a_1,\ldots,a_m) \in A^{m} ~\mid~ \delta_{A,A_{(a_3-a_2,\ldots,a_m-a_2)}}(a_1-a_2)\ge \tau \Bigr\} \Bigr| \,.
$$
Since $A$ is of SzT-type, we have, using \eqref{f:SzT-type}, 
\begin{align}\label{aa}
\sigma
&=\sum_{a_2,\dots,a_m\in A}\Bigl|\Bigl\{ a \in A ~\mid~ \delta_{A,A_{(a_3-a_2,\ldots,a_m-a_2)}}(a-a_2)\ge \tau \Big\} \Big|
\nonumber\\
&\le \sum_{a_2,\dots,a_m\in A}\Bigl|\Bigl\{ s \in A-A ~\mid~ \delta_{A,A_{(a_3-a_2,\ldots,a_m-a_2)}}(s)\ge \tau \Big\} \Big|
\nonumber\\
&\le
\frac{d(A) |A|}{\tau^3} \sum_{a_2,\dots,a_m\in A} |A_{(a_3-a_2,\ldots,a_m-a_2)}|^2 \\
& \le
\frac{d(A) |A|  \E_{m-1,3}(A)}{\tau^3}\,.\nonumber
\end{align}
Indeed, to see the last inequality, note that
the quantity $\sum_{a_2,\dots,a_m\in A}|A_{(a_3-a_2,\ldots,a_m-a_2)}|^2$ can be interpreted as the number of 
$(m+1)$-tuples $(a_2,\dots,a_m,z,z')$, with $a_2,\ldots,a_m\in A$ and
$z,z'\in A_{(a_3-a_2,\ldots,a_m-a_2)}$. 
We want to upper bound this quantity by $\E_{m-1,3}(A)$, which, by definition, is equal to 
$\sum_{y_3,\ldots,y_m\in\R}|A_{(y_3,\ldots,y_m)}|^3$.
This latter quantity is equal to the number of $(m+1)$-tuples $(y_3,\ldots,y_m,a_2,z,z')$
such that $a_2,z,z'\in A_{(y_3,\ldots,y_m)}$.

With each tuple $\tau=(y_3,\ldots,y_m,a_2,z,z')$ of the latter kind
we associate the tuple $\tau^*=(a_2,a_2+y_3,\ldots,a_2+y_m,z,z')$.
Since $a_2\in A_{(y_3,\ldots,y_m)}$, we have that each of $a_2, a_3:=a_2+y_3,\ldots, a_m:=a_2+y_m$
is in $A$.
Similarly, since $z,z'\in A_{(y_3,\ldots,y_m)}$, we have, by definition $z,z'\in A_{(a_3-a_2,\ldots,a_m-a_2)}$.
That is $\tau^*$ is a tuple that contributes to the quantity in \eqref{aa}.
The converse implication, passing from a tuple $\tau^*$ to that is counted in \eqref{aa}
to a tuple that is counted in $\E_{m-1,3}(A)$ is obtained by reversing the argument just given, in a straightforward manner.
The desired inequality thus follows, and this completes the proof of the lemma.   $\hfill\Box$

\section{Proof of Theorem~\ref{SzT2}}\label{sec:prf}

In what follows we assume, as in the theorem, that $d(A)\le d(B)$.

With each point $p=(a, b,c)\in A^3$ we associate a plane $h(p)$ in $\R^3$, given by
\begin{equation}\label{det}
\frac12\left|
\begin{array}{ccc}
a& b& c\\
x& y& z\\
1&1&1
\end{array}\right|=1,
\end{equation}
or equivalently by
\begin{equation*}
(c-b)x+(a-c)y+(b-a)z=2.
\end{equation*}
We put
$$
H:=\{h(p)\;\mid\;p\in A^3\}.
$$

A triangle with vertices $(a_1,x_1),(a_2,x_2),(a_3,x_3)$ has unit area if and only if the left-hand side of (\ref{det}) has absolute value 1,
so for half of the permutations $(i_1,i_2,i_3)$ (i.e., three permutations) of $(1,2,3)$, we have $(x_{i_1},x_{i_2},x_{i_3})\in h(a_{i_1},a_{i_2},a_{i_3})$.
In other words, the number of unit-area triangles is at most one third of the number of incidences between the points of $B^3$ and the planes of $H$.
In addition to the usual problematic issue that arises in point-plane incidence problems,
where many planes can pass through a line that contains many points (see, e.g., \cite{AS05}),
we need to face here the issue that the planes of $H$ are in general not distinct, and may arise with large multiplicity. Denote by $w(h)$ the \emph{multiplicity} of a plane $h\in H$, that is, $w(h)$ is the number of points $p\in A^3$ for which $h(p)=h$. Observe that, for $p,p'\in A^3$,
\begin{equation}
h(p)\equiv h(p')~\text{ if and only if }~p'\in p+(1,1,1)\R.
\end{equation}

We can transport this notion to points of $A^3$, by defining the \emph{multiplicity} $w(p)$ of a point $p\in A^3$ by
$$
w(p):= \Bigl|(p+(1,1,1)\R)\cap A^3\Bigr|.
$$
Then we clearly have $w(h(p))=w(p)$ for each $p\in A^3$.
Similarly, for $q\in B^3$, we put, by a slight abuse of notation,
$$
w(q):= \Bigl|(q+(1,1,1)\R)\cap B^3\Bigr|,
$$
and refer to it as the \emph{multiplicity} of $q$. (Clearly, the points of $B^3$ are all distinct, but the notion of their ``multiplicity" will become handy in one of the steps of the analysis --- see below.)

Fix a parameter $k\in\N$, whose specific value will be chosen later. We say that $h\in H$ (resp., $p\in A^3$, $q\in B^3$) is $k$-\emph{rich}, if its multiplicity is at least $k$; otherwise we say that it is $k$-\emph{poor}. For a unit-area triangle $T$, with vertices $(a,x), (b,y), (c,z)$, we say that $T$ is \emph{rich-rich} (resp., \emph{rich-poor}, \emph{poor-rich}, \emph{poor-poor}) if $(a,b,c)\in A^3$ is $k$-rich (resp., rich, poor, poor), and $(x,y,z)\in B^3$ is $k$-rich (resp., poor, rich, poor).
(These notions depend on the parameter $k$, which is fixed throughout the rest of the analysis.)

Next, we show that our assumption that $A$ and $B$ are convex, or, more generally, of SzT-type (with $\alpha=2$),
allows us to have some control on the multiplicity of the points and the planes, which we need for the proof.

%%%%%%%%%%%%%%%%%%%%%%%%%%%%%%%%%%%%%%%%%%%%%%%%%%%%%%
\begin{lemma}\label{krich}
Let $A$ be a set of $SzT$-type with parameter $\alpha=2$, and assume $|A|=n^{1/2}$.
Then the number of $k$-rich points in $A^3$ is $\displaystyle O\left(\frac{d^2 (A)n^2\log n}{k^3}\right)$.
\end{lemma}
%%%%%%%%%%%%%%%%%%%%%%%%%%%%%%%%%%%%%%%%%%%%%%%%%%%%%%
\noindent{\bf Proof.}
Note that the number of $k$-rich points in $A^3$ is exactly
$$
\Bigl|\left\{(a,b,c)\in A^3\mid \bigl|(A-a)\cap(A-b)\cap (A-c)\bigr|\ge k\right\}\Bigr|.
$$
Applying Lemma~\ref{mtaurich}, with the parameters $m=3$ and $\tau=k$, and then Lemma~\ref{E3}, the assertion follows.
$\hfill\Box$

We will also need the following lemma.
%%%%%%%%%%%%%%%%%%%%%%%%%%%%%%%%%%%%%%%%%%%%%%%%%%%%%%
\begin{lemma}\label{kproj}
Let $A$ be a set of $SzT$-type, with the parameter $\alpha=2$, and assume $|A|=n^{1/2}$.
The image set of the $k$-rich points of $A^3$, under the projection map onto the $xy$-plane, is of cardinality $O(d(A){n^{3/2}}/{k^2})$.
\end{lemma}
%%%%%%%%%%%%%%%%%%%%%%%%%%%%%%%%%%%%%%%%%%%%%%%%%%%%%%
\noindent{\bf Proof.}
First note that, by the definition of
SzT-type, the number of points $(a,b)\in A^2$,
for which the line $(a,b)+(1,1)\R$ contains at least $k$ points of $A^2$, is $O(d(A){n^{3/2}}/{k^2})$.
Indeed, the number of differences $s\in A-A$ with $\delta_{A,A}(s)\ge\tau$ is at most $d(A)n^{3/2}/\tau^3$.
Each difference $s$ determines, in a 1-1 manner,  a line in $\R^2$ with orientation $(1,1)$ that
contains the $\delta_{A,A}(s)$ pairs $(a,b)\in A^2$ with $b-a=s$.
Let $M_\tau$ (resp., $M_{\ge \tau}$) denote the number of differences $s\in A-A$
with $\delta_{A,A}(s)=\tau$ (resp., $\delta_{A,A}(s)\ge \tau$). Then the desired number of
points is
$$
\sum_{\tau\ge k}\tau M_\tau
=kM_{\ge k}+\sum_{\tau>k}M_{\ge \tau}
=O(d(A)n^{3/2}/k^2)+\sum_{\tau>k}O(d(A)n^{3/2}/\tau^3)
=O(d(A)n^{3/2}/k^2).
$$
Let $(a,b,c)\in A^3$ be $k$-rich. Then, by definition, the line $l:=(a,b,c)+(1,1,1)\R$ contains at least $k$ points of $A^3$.
We consider the line $l':=(a,b)+(1,1)\R$, which is the projection of $l$ onto the $xy$-plane, which we identify with $\R^2$.
Note that the projection of the points of $l\cap A^3$ onto $\R^2$ is injective and its image is contained in $l'\cap A^2$.
In particular, $l'$ contains at least $k$ points of $A^2$. As just argued, the total number of such points in $A^2$
(lying on some line of the form $l'$, that contains at least $k$ points of $A^2$) is $O(d(A){n^{3/2}}/{k^2})$.
$\hfill\Box$

In what follows, we bound separately the number of unit-area  triangles that are rich-rich, poor-rich (and, symmetrically, rich-poor), and poor-poor.

%RICH-RICH
\paragraph{Rich-rich triangles.}

Note that for $((a,b,c),(\xi,\eta))\in A^3\times B^2$, with $a\neq b$, there exists at most one point $\zeta\in B$
such that $T((a,\xi),(b,\eta),(c,\zeta))$ has unit area. Indeed, the point $(c,\zeta)$ must lie on a certain line
$l((a,\xi),(b,\eta))$ parallel to $(a,\xi)-(b,\eta)$. This line intersects $x=c$ in exactly one point (because
$a\neq b$), which determines the potential value of $\zeta$. Thus, since we are now concerned with the number of
rich-rich triangles (and focusing at the moment on the case where $a\neq b$), it suffices to bound the number of
such pairs $((a,b,c),(\xi,\eta))$, with $(a,b,c)\in A^3$ being rich, and $(\xi,\eta)\in B^2$ being the projection
of a rich point of $B^3$, which is
$$
O\left(\frac{d^2(A)n^2\log n}{k^3}\cdot \frac{d(B)n^{3/2}}{k^2}\right)=O\left( \frac{d^2(A)d(B)n^{7/2}\log n}{k^5}\right),
$$
using Lemma~\ref{krich} and Lemma~\ref{kproj}.

It is easy to check that the number of unit-area triangles $T(p,q,r)$, where $p,q,r\in P$ and $p,q$ share the
same abscissa (i.e., $A$-component), is $O(n^2)$. Indeed, there are $\Theta(n^{3/2})$ such pairs $(p,q)$, and
for each of them there exist at most $n^{1/2}$ points $r\in P$, such that $T(p,q,r)$ has unit area (because the
third vertex $r$ must lie on a certain line $l(p,q)$, which passes through at most this number of points of the
grid $P$); here we do not use the fact that we are interested only in rich-rich triangles.

We thus obtain the following lemma.
\begin{lemma}\label{rr}
The number of rich-rich triangles spanned by $P$ is $\displaystyle O\left( \frac{d^2(A)d(B)n^{7/2}\log n}{k^5} + n^2\right)$.
\end{lemma}

%POOR-RICH RICH-POOR
\paragraph{Poor-rich and rich-poor triangles.}
Consider first the case of poor-rich triangles. Put
$$
H_i:=\{h\in H\;\mid\;2^{i-1}\le w(h)< 2^i\},
$$
for $i=1,\ldots,\log k$, and
$$
S_{\ge k}:=\{q\in B^3\;\mid\;w(q)\ge k\}.
$$
That is, by definition, $\bigcup_iH_i$ is the collection of $k$-poor planes of $H$, and $S_{\ge k}$ is the set of $k$-rich points of $B^3$.
Since each element of $H_i$ has multiplicity at least $2^{i-1}$, we have the trivial bound $|H_i|\le n^{3/2}/2^{i-1}$.

Consider the family of horizontal planes ${\cal F}:=\{\xi_z\}_{z\in B}$, where $\xi_{z_0}:=\{z=z_0\}$.
Our strategy is to restrict $S_{\ge k}$ and $H_i$, for each fixed $i=1,\ldots,\log k$, to the planes $\xi\in{\cal F}$, and apply
the Szemer\'edi--Trotter incidence bound to the resulting collections of points and intersection lines, on each such $\xi$.
Note that two distinct planes $h_1,h_2\in H$ restricted to $\xi$, become two distinct lines in $\xi$. Indeed, each plane of $H$ contains a line
parallel to $(1,1,1)$, and two such planes, that additionally share a horizontal line within $\xi$, must be identical.
By definition, we have that $S_{\ge k}=\bigcup_{z\in B} (S_{\ge k}\cap \xi_z)$, and the union is of pairwise disjoint sets.
In particular, using Lemma~\ref{krich}, we have
$\displaystyle \sum_{z\in B}|S_{\ge k}\cap \xi_z|=|S_{\ge k}|=O\left(\frac{d^2(B)n^2\log n}{k^{3}}\right)$.

The number of incidences between the points of $S_{\ge k}$ and the poor planes of $H$, counted with multiplicity (of the planes), is at most
$$
\sum_{z\in B}\sum_{i=1}^{\log k}2^i\cdot {\cal{I}}(S_{\ge k}\cap \xi_z,H_{iz}),
$$
where $H_{iz}$ is the collection of lines $\{h\cap \xi_z\mid h\in H_i\}$. By Theorem~\ref{ST}, this is at most
\begin{align*}
\sum_{z\in B}\sum_{i=1}^{\log k}
2^i&\cdot O\left( |S_{\ge k}\cap \xi_z|^{2/3}\left(\frac{n^{3/2}}{2^{i-1}}\right)^{2/3}+ |S_{\ge k}\cap \xi_z|+\frac{n^{3/2}}{2^{i-1}}\right)\\
&=
O\left(
n \sum_{z\in B}|S_{\ge k}\cap \xi_z|^{2/3}\sum_{i=1}^{\log k}2^{i/3}
+ \sum_{z\in B}|S_{\ge k}\cap \xi_z|\sum_{i=1}^{\log k}2^i
+n^2\log k\right)\\
&=
O\left(
nk^{1/3} \sum_{z\in B}|S_{\ge k}\cap \xi_z|^{2/3}
+ \frac{d^2(B)n^2\log n}{k^2}
+n^2\log k\right)\\
&=
O\left(
nk^{1/3} \left(\sum_{z\in B}|S_{\ge k}\cap \xi_z|\right)^{2/3}|B|^{1/3}
+ \frac{d^2(B)n^2\log n}{k^2}
+n^2\log k\right)\\
&=
O\left(
\frac{d(B)^{4/3}n^{5/2}\log^{2/3}n}{k^{5/3}}
+ \frac{d^2(B)n^2\log n}{k^2}
+n^2\log k\right) ,
\end{align*}
where the second-to-last inequality is obtained via H\"older's inequality.

This bounds the number of poor-rich triangles spanned by $P$. Clearly, using a symmetric argument, in which
the roles of $A$ and $B$ are switched, this bound, with $d(A)$ replacing $d(B)$, also applies to the number
of rich-poor triangles spanned by $P$. We thus obtain the following lemma.
%%%%%%%%%%%%%%%%%%%%%%%%%%%%%%%%%%%%%%%%%%%%%%%
\begin{lemma}\label{pr}
The number of poor-rich triangles and of rich-poor triangles spanned by $P$ is
$\displaystyle O\left( \frac{(d(A)^{4/3}+d(B)^{4/3})n^{5/2}\log^{2/3}n}{k^{5/3}}+
\frac{(d^2(A)+d^2(B))n^2\log n}{k^2} + n^2\log k \right)$.
\end{lemma}
%%%%%%%%%%%%%%%%%%%%%%%%%%%%%%%%%%%%%%%%%%%%%%%
%%%%%%%%% POOR-POOR
\paragraph{Poor-poor triangles.}
Again we are going to use Theorem \ref{ST}.  For $i=1,\ldots,\log k$, put
$$
S_i:=\{q\in B^3\;\mid\;2^{i-1}\le w(q)< 2^i\},
$$
and let $S_i'$, $H_i'$ be the respective (orthogonal) projections of $S_i$, $H_i$ to the plane $\eta:=\{x+y+z=1\}$.
Note that $H_i'$ is a collection of lines in $\eta$.
Moreover, arguing as above, two distinct planes of $H_i$ project to two distinct lines of $H_i'$, and thus the multiplicity
of the lines is the same as the multiplicity of the original planes of $H_i$. Similarly, a point $q\in S_i$ with multiplicity
$t$ projects to a point $q'\in S_i'$ with multiplicity $t$ (by construction, there are exactly $t$ points of $S_i$ that project
to $q'$).
These observations allow us to use here, as before, the trivial bounds $|S_i'|\le n^{3/2}/2^{i-1}$, $|H_i'|\le n^{3/2}/2^{i-1}$,
for $i=1,\ldots,\log k$.

Applying Theorem~\ref{ST} to the collections $S_i',H_j'$ in $\eta$, for $i,j=1,\ldots,\log k$, taking into account the
multiplicity of the points and of the lines in these collections, we obtain that the number of incidences between the poor
points and the poor planes, counted with the appropriate multiplicity, is at most
\begin{align*}
\sum_{i,j=1}^{\log k}2^{i+j}\cdot I(S_i',H_j')
&=
\sum_{i,j=1}^{\log k}2^{i+j}\cdot
O\left(
\left(\frac{n^{3/2}}{2^{i-1}}\right)^{2/3}\left(\frac{n^{3/2}}{2^{j-1}}\right)^{2/3}+\frac{n^{3/2}}{2^{i-1}}+\frac{n^{3/2}}{2^{j-1}}
\right)\\
&=O\left(
n^2\sum_{i,j=1}^{\log k} 2^{(i+j)/3}+n^{3/2}\sum_{i,j=1}^{\log k}\left(2^i+2^j\right)
\right)
\\
&=
O\left(
n^2k^{2/3}+n^{3/2}k\log k
\right).
\end{align*}
Thus, we obtain the following lemma.
\begin{lemma}\label{pp}
The number of poor-poor triangles spanned by $P$ is $O\left(n^2k^{2/3}+n^{3/2}k\log k \right)$.
\end{lemma}

In summary, the number of unit-area triangles spanned by $P$ is
\begin{align}
O\Bigg(\frac{d^2(A)d(B)n^{7/2}}{k^5}\log n
&+\frac{(d(A)^{4/3}+d(B)^{4/3})n^{5/2}}{k^{5/3}}\log^{2/3}n\nonumber \\
&+\frac{(d^2(A)+d^2(B))n^2\log n}{k^2}
+n^2k^{2/3}+n^{3/2}k\log k\Bigg).
%\,.
\end{align}
Setting $k=n^{9/34}(d^2(A)d(B)\log n)^{3/17}$ and recalling that $d(A) \le d(B) \le n^{1/2}$ as well as our assumption 
$d(B) \le n^{6/47} d^{42/47} (A) (\log n)^{-13/47}$ makes the bound 
$\displaystyle O\left(n^{37/17}(d^2(A)d(B)\log n)^{2/17}\right)$, and Theorem \ref{SzT2} follows.
$\hfill\Box$

\end{document}